\magnification=1200

\catcode`\À=\active \defÀ{\`A}    \catcode`\à=\active \defà{\`a} 
\catcode`\Â=\active \defÂ{\^A}    \catcode`\â=\active \defâ{\^a} 
\catcode`\Æ=\active \defÆ{\AE}    \catcode`\æ=\active \defæ{\ae}
\catcode`\Ç=\active \defÇ{\c C}   \catcode`\ç=\active \defç{\c c}
\catcode`\È=\active \defÈ{\`E}    \catcode`\è=\active \defè{\`e} 
\catcode`\É=\active \defÉ{\'E}    \catcode`\é=\active \defé{\'e} 
\catcode`\Ê=\active \defÊ{\^E}    \catcode`\ê=\active \defê{\^e} 
\catcode`\Ë=\active \defË{\"E}    \catcode`\ë=\active \defë{\"e} 
\catcode`\Î=\active \defÎ{\^I}    \catcode`\î=\active \defî{\^\i}
\catcode`\Ï=\active \defÏ{\"I}    \catcode`\ï=\active \defï{\"\i}
\catcode`\Ô=\active \defÔ{\^O}    \catcode`\ô=\active \defô{\^o} 
\catcode`\Ù=\active \defÙ{\`U}    \catcode`\ù=\active \defù{\`u} 
\catcode`\Û=\active \defÛ{\^U}    \catcode`\û=\active \defû{\^u} 
\catcode`\Ü=\active \defÜ{\"U}    \catcode`\ü=\active \defü{\"u} 

\catcode`\ =\active \def { }

\hsize=11.25cm    
\vsize=18cm       
\parindent=12pt   \parskip=5pt     

\hoffset=.5cm   
\voffset=.8cm   

\pretolerance=500 \tolerance=1000  \brokenpenalty=5000

\catcode`\@=11

\font\eightrm=cmr8         \font\eighti=cmmi8
\font\eightsy=cmsy8        \font\eightbf=cmbx8
\font\eighttt=cmtt8        \font\eightit=cmti8
\font\eightsl=cmsl8        \font\sixrm=cmr6
\font\sixi=cmmi6           \font\sixsy=cmsy6
\font\sixbf=cmbx6

\font\tengoth=eufm10 
\font\eightgoth=eufm8  
\font\sevengoth=eufm7      
\font\sixgoth=eufm6        \font\fivegoth=eufm5

\skewchar\eighti='177 \skewchar\sixi='177
\skewchar\eightsy='60 \skewchar\sixsy='60

\newfam\gothfam           \newfam\bboardfam

\def\tenpoint{
  \textfont0=\tenrm \scriptfont0=\sevenrm \scriptscriptfont0=\fiverm
  \def\rm{\fam\z@\tenrm}
  \textfont1=\teni  \scriptfont1=\seveni  \scriptscriptfont1=\fivei
  \def\oldstyle{\fam\@ne\teni}\let\old=\oldstyle
  \textfont2=\tensy \scriptfont2=\sevensy \scriptscriptfont2=\fivesy
  \textfont\gothfam=\tengoth \scriptfont\gothfam=\sevengoth
  \scriptscriptfont\gothfam=\fivegoth
  \def\goth{\fam\gothfam\tengoth}
  
  \textfont\itfam=\tenit
  \def\it{\fam\itfam\tenit}
  \textfont\slfam=\tensl
  \def\sl{\fam\slfam\tensl}
  \textfont\bffam=\tenbf \scriptfont\bffam=\sevenbf
  \scriptscriptfont\bffam=\fivebf
  \def\bf{\fam\bffam\tenbf}
  \textfont\ttfam=\tentt
  \def\tt{\fam\ttfam\tentt}
  \abovedisplayskip=12pt plus 3pt minus 9pt
  \belowdisplayskip=\abovedisplayskip
  \abovedisplayshortskip=0pt plus 3pt
  \belowdisplayshortskip=4pt plus 3pt 
  \smallskipamount=3pt plus 1pt minus 1pt
  \medskipamount=6pt plus 2pt minus 2pt
  \bigskipamount=12pt plus 4pt minus 4pt
  \normalbaselineskip=12pt
  \setbox\strutbox=\hbox{\vrule height8.5pt depth3.5pt width0pt}
  \let\bigf@nt=\tenrm       \let\smallf@nt=\sevenrm
  \normalbaselines\rm}

\def\eightpoint{
  \textfont0=\eightrm \scriptfont0=\sixrm \scriptscriptfont0=\fiverm
  \def\rm{\fam\z@\eightrm}
  \textfont1=\eighti  \scriptfont1=\sixi  \scriptscriptfont1=\fivei
  \def\oldstyle{\fam\@ne\eighti}\let\old=\oldstyle
  \textfont2=\eightsy \scriptfont2=\sixsy \scriptscriptfont2=\fivesy
  \textfont\gothfam=\eightgoth \scriptfont\gothfam=\sixgoth
  \scriptscriptfont\gothfam=\fivegoth
  \def\goth{\fam\gothfam\eightgoth}
  
  \textfont\itfam=\eightit
  \def\it{\fam\itfam\eightit}
  \textfont\slfam=\eightsl
  \def\sl{\fam\slfam\eightsl}
  \textfont\bffam=\eightbf \scriptfont\bffam=\sixbf
  \scriptscriptfont\bffam=\fivebf
  \def\bf{\fam\bffam\eightbf}
  \textfont\ttfam=\eighttt
  \def\tt{\fam\ttfam\eighttt}
  \abovedisplayskip=9pt plus 3pt minus 9pt
  \belowdisplayskip=\abovedisplayskip
  \abovedisplayshortskip=0pt plus 3pt
  \belowdisplayshortskip=3pt plus 3pt 
  \smallskipamount=2pt plus 1pt minus 1pt
  \medskipamount=4pt plus 2pt minus 1pt
  \bigskipamount=9pt plus 3pt minus 3pt
  \normalbaselineskip=9pt
  \setbox\strutbox=\hbox{\vrule height7pt depth2pt width0pt}
  \let\bigf@nt=\eightrm     \let\smallf@nt=\sixrm
  \normalbaselines\rm}

\tenpoint

\def\pc#1{\bigf@nt#1\smallf@nt}         \def\pd#1 {{\pc#1} }

\catcode`\;=\active
\def;{\relax\ifhmode\ifdim\lastskip>\z@\unskip\fi
\kern\fontdimen2  -1.2 \fontdimen3 \string;}

\catcode`\:=\active
\def:{\relax\ifhmode\ifdim\lastskip>\z@\unskip\fi\penalty\@M\ \fi\string:}

\catcode`\!=\active
\def!{\relax\ifhmode\ifdim\lastskip>\z@
\unskip\fi\kern\fontdimen2  -1.1 \fontdimen3 \string!}

\catcode`\?=\active
\def?{\relax\ifhmode\ifdim\lastskip>\z@
\unskip\fi\kern\fontdimen2  -1.1 \fontdimen3 \string?}

\catcode`\«=\active 
\def«{\raise.4ex\hbox{%
 $\scriptscriptstyle\langle\!\langle$}}

\catcode`\»=\active 
\def»{\raise.4ex\hbox{%
 $\scriptscriptstyle\rangle\!\rangle$}}

\frenchspacing

\def\raggedbottom{\topskip 10pt plus 36pt\r@ggedbottomtrue}

\def\pointir{\unskip . --- \ignorespaces}

\def\Medbreak{\vskip-\lastskip\medbreak}

\long\def\th#1 #2\enonce#3\endth{
   \Medbreak\noindent
   {\pc#1} {#2\unskip}\pointir{\it #3}\smallskip}

\def\proof{\vskip-\lastskip\smallskip\noindent
 {\it Proof} : }

\def\decale#1{\smallbreak\hskip 28pt\llap{#1}\kern 5pt}
\def\decaledecale#1{\smallbreak\hskip 34pt\llap{#1}\kern 5pt}
\def\puce{\smallbreak\hskip 6pt{$\scriptstyle\bullet$}\kern 5pt}

\def\eqalign#1{\null\,\vcenter{\openup\jot\m@th\ialign{
\strut\hfil$\displaystyle{##}$&$\displaystyle{{}##}$\hfil
&&\quad\strut\hfil$\displaystyle{##}$&$\displaystyle{{}##}$\hfil
\crcr#1\crcr}}\,}

\catcode`\@=12

\showboxbreadth=-1  \showboxdepth=-1

\newcount\numerodesection \numerodesection=1
\def\section#1{\bigbreak
 {\bf\number\numerodesection.\ \ #1}\nobreak\medskip
 \advance\numerodesection by1}

\mathcode`A="7041 \mathcode`B="7042 \mathcode`C="7043 \mathcode`D="7044
\mathcode`E="7045 \mathcode`F="7046 \mathcode`G="7047 \mathcode`H="7048
\mathcode`I="7049 \mathcode`J="704A \mathcode`K="704B \mathcode`L="704C
\mathcode`M="704D \mathcode`N="704E \mathcode`O="704F \mathcode`P="7050
\mathcode`Q="7051 \mathcode`R="7052 \mathcode`S="7053 \mathcode`T="7054
\mathcode`U="7055 \mathcode`V="7056 \mathcode`W="7057 \mathcode`X="7058
\mathcode`Y="7059 \mathcode`Z="705A


\catcode`\À=\active \defÀ{\`A}    \catcode`\à=\active \defà{\`a} 
\catcode`\Â=\active \defÂ{\^A}    \catcode`\â=\active \defâ{\^a} 
\catcode`\Æ=\active \defÆ{\AE}    \catcode`\æ=\active \defæ{\ae}
\catcode`\Ç=\active \defÇ{\c C}   \catcode`\ç=\active \defç{\c c}
\catcode`\È=\active \defÈ{\`E}    \catcode`\è=\active \defè{\`e} 
\catcode`\É=\active \defÉ{\'E}    \catcode`\é=\active \defé{\'e} 
\catcode`\Ê=\active \defÊ{\^E}    \catcode`\ê=\active \defê{\^e} 
\catcode`\Ë=\active \defË{\"E}    \catcode`\ë=\active \defë{\"e} 
\catcode`\Î=\active \defÎ{\^I}    \catcode`\î=\active \defî{\^\i}
\catcode`\Ï=\active \defÏ{\"I}    \catcode`\ï=\active \defï{\"\i}
\catcode`\Ô=\active \defÔ{\^O}    \catcode`\ô=\active \defô{\^o} 
\catcode`\Ù=\active \defÙ{\`U}    \catcode`\ù=\active \defù{\`u} 
\catcode`\Û=\active \defÛ{\^U}    \catcode`\û=\active \defû{\^u} 
\catcode`\Ü=\active \defÜ{\"U}    \catcode`\ü=\active \defü{\"u}

\centerline{\bf The Tao of Mathematics}
\bigskip
\centerline{\it Notes for a lecture at the}
\smallskip
\centerline{\it Indian Statistical Institute, Bangalore, 9 August,
2005.}
\bigskip

{\bf 1.  Ramanujan's $\tau$-function.}  I've come across your syllabus
and I was impressed by the fact that Ramanujan's $\tau$-function is
taught to undergraduates.  His famous work on this function dates back
to 1916.  As you know, the function $\tau:{\bf
N}^\times\rightarrow{\bf Z}$ is defined formally as the coefficients
in the power-series expansion of
$$\eqalign{
\Delta(q)&=q\prod_{n\ge1}(1-q^n)^{24}\cr
 &=\sum_{n\ge1}\tau(n)q^n,\cr}
$$
where $q$ is an indeterminate.  Ramanujan made three conjectures about
$\tau$:
$$
\cases{
\tau(mn)=\tau(m)\tau(n)\quad\hbox{if pgcd}(m,n)=1,\cr
\tau(p^{\alpha+2})=\tau(p)\tau(p^{\alpha+1})-p^{11}\tau(p^\alpha)
\quad(p \hbox{ prime},\ \alpha\ge0)\cr}
\leqno{\hbox{Conjecture I}}
$$\vskip-10pt
$$ |\tau(p)|\leq 2p^{11\over2}\quad(p \hbox{ prime}).
\leqno{\hbox{Conjecture II}}
$$
This is the same as saying that the reciprocals $\alpha$, $\beta$ of
the roots of the polynomial $1-\tau(p)T+p^{11}T^2$ satisfy
$|\alpha|=|\beta|=p^{11\over2}$.

Then there were many congruences (modulo $2^{11}$, $3^7$, $5^3$, $7$,
$23$ and $691$), some of which he proved, for example
$$
\tau(p)\equiv1+p^{11}\pmod{691}\quad\hbox{($p$ prime $\neq691$)}.
\leqno{\hbox{Conjecture III}}
$$
It seems that Hardy, when talking about this paper of Ramanujan, says
that it belongs to the backwaters of mathematics.  Weil did not fail
to point out that this observation of Hardy merely shows the
difference in taste between analysts and arithmeticians.

Conjecture I was proved by Mordell soon afterwards (1918).  It is now
a consequence of a general theory developed by Hecke.  Indeed,
Mordell's proof can be viewed as a precocious use of Hecke operators.

Conjecture II became part of a general theory when it was observed by
Ihara that it follows from Weil's conjectures on the zeta functions of
varieties over finite fields, which were finally proved by Deligne in
1973.  It holds the world record for the ratio
$$
{\hbox{Length of the proof}\over
\hbox{Length of the statement}}.
$$

Serre was prompted by Ramanujan's congruences to come up with his
conjectures relating modular forms and representations of ${\rm
Gal}(\bar{\bf Q}|{\bf Q})$; this is going to be our main concern
today.

Many of the congruences in Ramanujan's conjecture III were proved by
Bambah ($\tau(p)\equiv1+p^{11}\pmod{2^5}$, $p\neq2$), K.~G.~Ramanathan
($\tau(p)\equiv1+p\pmod{3}$, $p\neq3$), and others.

A systematic theory was developed by Serre and Swinnerton-Dyer.  They
constructed, for each prime $l$, a representation $\rho_l:{\rm Gal}(\bar{\bf
Q}|{\bf Q})\rightarrow {\rm GL}_2({\bf F}_l)$ which is unramified
outside $l$ and such that for every prime $p\neq l$, one has
$$
{\rm Tr}(\rho_l({\rm Frob}_p))\equiv\tau(p)
\quad{\rm and}\quad
{\rm det}(\rho_l({\rm Frob}_p))\equiv p^{11}\quad(\hbox{in }{\bf F}_l),
$$
where ${\rm Frob}_p$ is a lift of the automorphism $x\mapsto x^p$ of
${\bf F}_p$.  In this way, they were able to give a uniform proof of
all of Ramanujan's congruences.  They could also show that no further
congruences hold for the $\tau$-function.

Now it so happens that the more fundamental function is not $\tau$ but
$\Delta$.  Putting $q=e^{2i\pi z}$, we can view $\Delta$ as a function
of $z$ (in the upper half-plane ${\rm Im}(z)>0$); it has some amazing
properties.  For example, we know that ${\bf SL}_2({\bf Z})$ acts on
the upper half-plane in a natural manner :
$$
\gamma.z={az+b\over cz+d},\qquad
\gamma=\pmatrix{a&b\cr c&d\cr}\in{\bf SL}_2({\bf Z}).
$$
It can be verified that $\Delta(\gamma.z)=(cz+d)^{12}\Delta(z)$ for
every $\gamma$ and every $z$; this is expressed by saying that
$\Delta$ is a cusp form of weight $12$ for ${\bf SL}_2({\bf Z})$.

{\bf 2.  Modular forms.}  Quite generally, an analytic function $f$ on
the upper half-plane is said to be modular of weight $k$ for a
subgroup $\Gamma$ of ${\bf SL}_2({\bf Z})$ is one has
$$
f(\gamma.z)=(cz+d)^kf(z),\qquad\gamma=\pmatrix{a&b\cr c&d\cr}
$$
for every $z$ and every $\gamma\in\Gamma$.  If $f$ vanishes at the
``cusps'' of $\Gamma$, it is called a cusp form.  The most interesting
and useful case occurs when $\Gamma$ is a congruence subgroup of ${\bf
SL}_2({\bf Z})$, i.e.~defined by congruence conditions such as
$$
\pmatrix{a&b\cr c&d\cr}\equiv
\pmatrix{1&0\cr 0&1\cr}\pmod{N}
$$
for some fixed integer $N$, which is then called the ``level'' of
those modular forms.

Cusp forms of a given level~$N$ and weight~$k$ form a
finite-dimensional vector space over ${\bf C}$.  They have been
extensively studied ; tables are available online for example on
William Stein's webpage.  The space comes with a natural family
of commuting operators, called the Hecke operators, indexed by the
primes.  Cusp forms which are eigenvectors for all of these operators,
called eigenforms, are of special relevance ; they are the ones which
show up in so many places.  

For example, what Wiles and his followers really proved is that every
elliptic curve $E$ over ${\bf Q}$, i.e.~a curve of ``genus~1'' given
by an equation of the type
$$
y^2=x^3+ax+b\qquad (a,b\in{\bf Q}),
$$ 
is ``modular'', i.e.~there exists an eigenform $f$, of weight~$2$ and
appropriate level, whose eigenvalues $a_p$ satisfy : for almost every
prime $p$, there are precisely $p-a_p$ solutions of the congurence
$y^2\equiv x^3+ax+b\pmod{p}$.  Gauss was the first to count the number
of points on a specific elliptic curve modulo various primes ; Weil
the first to prove the modularity of a specific elliptic curve.

The simplest eigenform is $\Delta$ (level~$1$, weight~$k=12$), as was
known essentially to Jacobi :
$$
\Delta\left(\pmatrix{a&b\cr c&d\cr}z\right)=
(cz+d)^{12}\Delta(z),\qquad
\pmatrix{a&b\cr c&d\cr}\in{\bf SL}_2({\bf Z}).
$$
The eigenvalue $a_p$ is $\tau(p)$ for every prime~$p$.  This explains
why we have $11=k-1$ as exponent in Conjectures II and III.

We are lucky in that modular forms admit such a concrete description,
as analytic functions on the upper half-plane.  There is a more
abstract definition, in terms of representations of the ${\bf GL}_2$
of the ad{\`e}les of ${\bf Q}$, which makes it clear how they are the
correct generalisation in degree~$2$ of characters $({\bf Z}/N{\bf
Z})^\times\rightarrow{\bf C}^\times$ (degree~$1$), and what should
take their place in degree~$3$ etc.  

Another approach to modular forms, the geometric one, sees them as
sections of certain line bundles on ``modular curves'' ; it allows us
to define modular forms over arbitrary rings.  Indeed, what show up in
Serre's conjectures are modular forms over $\bar{\bf F}_l$.

{\bf 3.  Serre's conjectures.}  After the representation attached to
$\Delta$, Serre conjectured that there were similar representations
associated to every eigenform.  This was proved by Shimura in
weight~$2$ and by Deligne in general ; the weight~$1$ case is slightly
different : in place of representations into ${\bf GL}_2(\bar{\bf
Q}_l)$, we get (Deligne-Serre) representations into ${\bf GL}_2({\bf
C)}$.  Leaving aside this case, to every eigenform~$f$, of level~$N$
and weight~$k$, there corresponds a continuous representation
$\rho_f:{\rm Gal}(\bar{\bf Q}|{\bf Q})\rightarrow{\bf GL}_2(\bar{\bf
F}_l)$ which is unramified outside $Nl$ and such that for every prime
$p$ not dividing $Nl$, one has
$$
{\rm Tr}(\rho_f({\rm Frob}_p))\equiv a_p
\quad{\rm and}\quad
{\rm det}(\rho_f({\rm Frob}_p))\equiv p^{k-1}\quad(\hbox{in }\bar{\bf F}_l),
$$ 
where $f$ is sent to $a_pf$ by the $p^{\rm th}$ Hecke operator.
Having got this far, Serre tentatively talked about a converse.  

The converse says that every odd (meaning ${\rm det}(\rho(c))=-1$,
where $c\in{\rm Gal}(\bar{\bf Q}|{\bf Q})$ is a complex conjugation)
irreducible representation $\rho:{\rm Gal}(\bar{\bf Q}|{\bf
Q})\rightarrow{\bf GL}_2(\bar{\bf F}_l)$ arises from an eigenform by
the above construction.  He wrote about it to Tate in the early 70s ;
the reply was a proof for representations into ${\bf GL}_2({\bf
F}_2)$.  Using the same technique, Serre could prove the converse for
representations into ${\bf GL}_2({\bf F}_3)$.

Some time in the mid-80s, Colmez pointed out to Serre some of the
amazing consequences of his conjecture, among them the conjecture that
for every elliptic curve $E$ over ${\bf Q}$ is ``modular'', as
eventually proved by Wiles and his coworkers in full generality in the
late 90s.  These astonishing ({\it inqui\'etantes\/}) consequences of
his conjectures prompted Serre to make them more precise, by pinning
down the level and the weight of the eigenform giving rise to the
given odd irreducible representation.  What is nice about this
refinement is that it makes the conjecture verifiable.  Given $\rho$,
one computes the level $N(\rho)$ and the weight $k(\rho)$ and checks
in a finite amount of time if there is an eigenform of the type
required by the conjecture.

Results of Edixhoven, Ribet and others have shown that if there is
some eigenform giving rise to an odd irreducible $\rho:{\rm
Gal}(\bar{\bf Q}|{\bf Q})\rightarrow{\bf GL}_2(\bar{\bf F}_l)$, then
there is an eigenform $f$ of the correct weight $k(\rho)$ and the
correct level $N(\rho)$ giving rise to $\rho$.  So the refined version
of the conjecture follows from the na{\"\i}ve version.

{\bf 4.  Khare's proof in level 1.}  Following a stratagy he had
worked out at the end of last year with Wintenberger, the proof is an
elaborate induction on the prime~$l$.  The first few $l$ have to be
handled individually ; the induction works only for sufficiently
big~$l$.  It is an awsome achievement, using some of the deepest
results about galoisian representations obtained since the proof by
Wiles and others that every elliptic curve over ${\bf Q}$ is
``modular''.

{\bf 5.  Any questions~?}  You could ask : Why work with such a
complicated field as ${\bf Q}$, which has infinitely many primes ?
Why not first study things over a simpler field such as ${\bf Q}_p$,
which has only one prime ?  There are two local problems, the case
$l\neq p$ and the case $l=p$. Khare himself, and Marie-France
Vign{\'e}ras have proved the local analogue of Serre's conjecture when
$l\neq p$.  Progress has been made by Breuil on the case $l=p$, where
even the formulation of the problem is not easy.

You could also ask : why not first study (local and global)
representations into ${\bf GL}_1(\bar{\bf F}_1)=\bar{\bf
F}_l^\times$ ?  Indeed, this is what was done from Takagi to Artin ;
this is now a part of the theory of abelian extensions of local or
global fields.  Serre's conjecture is really the first step in the
programme (Langlands) of establishing a similar theory for all
galoisian extensions.
\goodbreak

{\bf 6.  A glimpse into the future.}  We began by recalling
Ramanujan's three conjectures about his $\Delta$.  We have seen how
all three are connected to some of the most profound mathematical
theories of the last century (Hecke theory, Weil conjectures, Langlands
programme).  

But there is more to Ramanujan's $\Delta$ ; in an article to be
written by the end of this year (2005), three authors (Edixhoven,
Couveignes and R.~de Jong) ``will clearly prove that the mod~$l$
Galois representations associated to the modular form $\Delta$ can be
computed in time polynomial in $l$''.

\bigskip

\vfill
\leftline{Chandan Singh Dalawat}
\leftline{Harish-Chandra Research Institute}
\leftline{Chhatnag Road, Jhunsi}
\leftline{Allahabad 211 019}
\vfill\eject

\centerline{\bf Think locally}
\bigskip
\centerline{\it Notes for a colloquium talk at the}
\smallskip
\centerline{\it Indian Statistical Institute, Bangalore, 9 August,
2005.}
\bigskip
\th THEOREM
\enonce
The polynomial\/ $x^2+1$ has no roots in\/ ${\bf Q}$.
\endth
\proof $x^2+1$ has no roots in ${\bf R}$.

\th THEOREM
\enonce
The polynomial\/ $x^2-2$ has no roots in\/ ${\bf Q}$.
\endth
The previous proof does not work, because $x^2-2$ does have (two)
roots in ${\bf R}$.  Any root in ${\bf Q}$ will have to belong to
${\bf Z}$.  For $\alpha\in{\bf Z}$, computing modulo~$5$ gives
$$\vbox{\halign{\hfil$#$\quad \ \ &\hfil$#$\quad\ \ &\hfil$#$\cr
\alpha&\alpha^2&\alpha^2-2\cr
\noalign{\vskip-5pt}
\multispan3\hrulefill\cr
0&0&-2\cr
1&1&-1\cr
2&-1&2\cr
-2&-1&2\cr
-1&1&-1\cr
\noalign{\vskip-5pt}
\multispan3\hrulefill\cr
}}$$ 
We see that $\alpha^2-2\not\equiv0\pmod{5}$.  Hensel
would have reformulated this little computation as :
\proof $x^2-2$ has no roots in ${\bf Q}_5$.   

In general, for every prime number $p$, we have a locally compact
field ${\bf Q}_p$ which plays a similar role : it can sometimes be
used to show that something doesn't happen over ${\bf Q}$ because it
doesn't happen over ${\bf Q}_p$.

The field ${\bf Q}_p$ can be defined in many equivalent ways.  The
simplest would be to say that it is the field of fractions of the ring
${\bf Z}_p$, which is integral.  One could also define it as the
completion of ${\bf Q}$ with respect to the distance $|x-y|_p$, where
the absolute value $|\phantom{x}|_p$ comes from the valuation $v_p:{\bf
Q}^\times\rightarrow{\bf Z}$ sending $x$ to the power of $p$ in $x$.

Hensel was pursuing an analogy between number fields and compact
connected analytic curves.  He viewed these valuations, or rather the
completions with respect to the associated absolute values, as being
the ``points'' on the ``curve'' that the number field ``is''.
Ostrowski showed that these are the only ``points'', i.e.~every
discrete valuation of ${\bf Q}$ comes from some prime number $p$.
This point of view is fully justified by Grothendieck's theory of
schemes, where there is indeed a ``curve'' (i.e.~a $1$-dimensional
scheme) corresponding to (the ring of integers of) a number field $K$
whose closed points are precisely the various discrete valuations of
$K$.  Arithmetic requires that the ``places at infinity'', namely the
completion ${\bf R}$ when $K={\bf Q}$, be treated on an equal footing ;
the full significance of this requirement was realised by Arakelov.

I want to emphasize that arithmetical problems should be first studied
locally, i.e.~one place at a time.  The trend among analysts and
topologists is to seek global results ; often the local versions are
trivial.  Not so in arithmetic.  We too seek global results ---
applicable to number fields --- but our local problems are hardly ever
trivial.

As an illustration of this kind of thinking, let me recount my
encounter with one your students, Anupam Kumar Singh, this morning.
Along with your colleague Maneesh Thakur, he is studying ``reality''
of an element $x\in G(k)$ in the group of rational points of a linear
algebraic group $G$ over a field $k$ ; an element in a group is said
to be ``real'' if it is conjugate to its inverse.  They have examples
of $x\in G({\bf Q})$ which are not ``real''.  My first question was to
ask if that $x$ is ``real'' over every place of ${\bf Q}$.  If there is a
place $v$ where $x$ is not ``real'', we have a {\it stronger\/}
statement : not only is $x$ not ``real'' over ${\bf Q}$, it is not
``real'' over ${\bf Q}_v$.  If $x$ happens to be ``real'' at every
place, we have an instance of the failure of a {\it local-to-global
principle}.  Many people would find such an example interesting and
would like to look for obstructions to account for it.

Today, we are going to work over a local field such as ${\bf Q}_p$.
Suppose you are given a smooth projective variety over $X_\eta$ over
${\bf Q}_p$, i.e.~a variety definable by a system of homogenous
polynomials
$$
f_1=0,\; f_2=0,\; \ldots,\; f_r=0\ ;\quad\quad
f_i\in{\bf Q}_p[T_0,\ldots,T_n]
$$
with a condition ensuring that there are no singularities.
Multiplying all the $f_i$ by a suitable power of $p$, we may assume
that they have coefficients in ${\bf Z}_p$.  Reducing the $f_i$ modulo
$p$, we get a system
$$
\bar f_1=0,\; \bar f_2=0,\; \ldots,\; \bar f_r=0\ ;\quad\quad
f_i\in{\bf F}_p[T_0,\ldots,T_n]
$$
defining a variety $X_s$ over ${\bf F}_p$.  The variety $X_s$ need not
be smooth.  If it is smooth for some choice of $f_i$ defining
$X_\eta$, we say that $X_\eta$ has {\it good reduction} ; it is said
to have {\it bad reduction\/} otherwise.
\vfill
\centerline{$\langle${\it The rest came from the talk at
Madras, at a more explicit level.}$\rangle$}
\vfill
\leftline{Chandan Singh Dalawat}
\leftline{Harish-Chandra Research Institute}
\leftline{Chhatnag Road, Jhunsi}
\leftline{Allahabad 211 019}
\vfill\eject
\bye